%% file: periodic.tex
\documentclass[12pt,twoside]{article}
\usepackage{amsmath,amsfonts,amssymb,amsbsy,a4,enumerate,epsfig,array,theorem,psfrag}

\newtheorem{theo}{Theorem}
\newtheorem{fact}{Fact}
\newtheorem{prop}{Proposition}[section]
\newtheorem{coro}[prop]{Corollary}

\newtheorem{lemma}[prop]{Lemma}

\newcommand{\expg}{\exp_{\fg}}

\newcommand{\cone}{\bowtie}

\newcommand{\id}{\textrm{id}}

\newcommand{\interior}{\operatorname{int}}
\newcommand{\tr}{\operatorname{tr}}

\newcommand{\Lev}{\textrm{Lev}}

\newcommand{\fg}{{\mathfrak g}}

\newcommand{\ZZ}{{\mathbb{Z}}}
\newcommand{\RR}{{\mathbb{R}}}
\newcommand{\CC}{{\mathbb{C}}}
\newcommand{\Ss}{{\mathbb{S}}}
\newcommand{\BB}{{\mathbb{B}}}

\newcommand{\HH}{{\mathbb{H}}}
\newcommand{\Bb}{{\cal{B}}}
\newcommand{\Cc}{{\cal{C}}}
\newcommand{\Dd}{{\cal{D}}}

\newcommand{\Hh}{{\cal{H}}}
\newcommand{\Kk}{{\cal{K}}}

\newcommand{\Vv}{{\cal{V}}}
\newcommand{\Ww}{{\cal{W}}}

\newcommand{\vv}{{\mathbf{v}}}

\newcommand{\beps}{{\boldsymbol\epsilon}}
\newcommand{\bone}{{\boldsymbol 1}}

\newcommand{\nobf}{\noindent\bf}

\def\qed{\unskip\nobreak\hfil\penalty50\hskip1.75em\null\nobreak\hfil
$\blacksquare$ {\parfillskip=0pt \finalhyphendemerits=0 \par}\goodbreak}
\def\qede{\unskip\nobreak
\hfil\penalty50\hskip1.75em\null\nobreak\hfil$\square$
{\parfillskip=0pt \finalhyphendemerits=0 \par}}

\begin{document}
\title{The geometry of the critical set of \\
nonlinear periodic Sturm-Liouville operators}
\author{Dan Burghelea, Nicolau C. Saldanha and Carlos Tomei}
\maketitle

\input{periodicI}

\input{periodicA1}
\input{periodicX}

\end{document}

%% file: periodicI.tex
\begin{abstract}
We study the critical set $\Cc$
of the nonlinear differential operator $F(u) = -u'' + f(u)$
defined on a Sobolev space of periodic functions $H^p(\Ss^1)$, $p \ge 1$.
Let $\RR^2_{xy} \subset \RR^3$ be the plane $z = 0$ and,
for $n > 0$, let $\cone_n$ be the cone
$x^2 + y^2 = \tan^2 z$, $|z - 2\pi n| < {\pi}/{2}$;
also set $\Sigma = \RR^2_{xy} \cup \bigcup_{n > 0} \cone_n$.
For a generic smooth nonlinearity $f: \RR \to \RR$ with surjective derivative,
we show that there is a diffeomorphism
between the pairs $(H^p(\Ss^1), \Cc)$ and $(\RR^3, \Sigma) \times \HH$
where $\HH$ is a real separable infinite dimensional Hilbert space.
\end{abstract}

\medbreak


\smallbreak

{\noindent\bf MSC-class:} 34B15; 34B24; 46T05.

\section{Introduction}
\label{section:intro}

A basic object in the study of
a smooth nonlinear operator $F: X \to Y$ between Banach spaces
is its \textit{critical set} $\Cc \subset X$,
the set of points $u \in X$ for which the derivative $DF(x)$
is not a (bounded) linear isomorphism between $X$ and $Y$.
For instance, knowledge of $C$ and $F(C)$ yields substantial
information about the number of solutions of the equation $F(x) = b$,
$b \in Y$ (\cite{MST0}, \cite{MST1}).
In this paper we consider a special but relevant example,
the nonlinear periodic Sturm-Liouville operator and 
determine the topology of the pair $(X, \Cc)$ in the generic case.

More precisely, for a smooth nonlinearity $f: \RR \to \RR$
and $p \ge 1$ denote by $F$ the differential operator
\begin{align*} F: H^p(\Ss^1) &\to H^{p-2}(\Ss^1). \\
u &\mapsto -u'' + f(u) \end{align*}
Here $\Ss^1 = \RR/2\pi\ZZ$ and
$H^p(\Ss^1)$ is the Sobolev space of periodic functions $u(t)$ with
square integrable $p$-th derivative.
Clearly, $DF(u): H^p(\Ss^1) \to H^{p-2}(\Ss^1)$
is $DF(u) v = -v'' + f'(u) v$,
a Fredholm operator of index 0.
Thus, the critical set $\Cc \subset H^p(\Ss^1)$ of $F$ is
\[ \Cc = \{ u \in H^p(\Ss^1) \;|\;
DF(u) \textrm{ has nontrivial kernel} \}. \]
In other words, $u \in \Cc$ if and only if the equation
\[ -v'' + f'(u) v = 0, \quad v(0) = v(2\pi), \quad v'(0) = v'(2\pi) \]
admits a nonzero solution.

This paper continues the project of obtaining geometric understanding
for some nonlinear operators $F: X \to Y$.
The starting point might be located in a fundamental
result of Ambrosetti and Prodi (\cite{AP}; see also \cite{AM}).
As interpreted by Berger, Church and Podolak (\cite{BC}, \cite{BP}),
it states that, for appropriate nonlinearities,
$u \mapsto -\Delta u + f(u)$
acting on functions satisfying Dirichlet boundary conditions is a global fold.
Subsequently, a number of different operators were shown to be
either global folds or global cusps
(see, among others, \cite{CT}, \cite{MST1}, \cite{Ruf}).
Ideally, the description of the critical set $\Cc$ should
include its stratification into different kinds of singularities.

In many examples, $\Cc$ and its strata were shown to be topologically trivial.
Differential operators related to nonlinear
Sturm-Liouville second order problems $F(u) =  - u'' + f(u)$
in finite intervals have very different critical sets
depending on the boundary conditions.
For generic nonlinearities $f$, Dirichlet boundary conditions
give rise to a critical set $\Cc$ which is
ambient diffeomorphic to a countable (possibly finite)
union of parallel hyperplanes (\cite{BT}, \cite{BST1}).
The situation is very different for periodic boundary conditions.

We now describe the generic nonlinearities $f$ for which
our main result holds.
A smooth function $g: \RR \to \RR$ is {\it nowhere flat} if
there exists a positive integer $r$ such that for all $x \in \RR$,
there exists $r' = r'(x) \in \ZZ$, $0 < r' \le r$, such that
$g^{(r')}(x) \ne 0$.
An \textit{admissible nonlinearity} is a smooth function
$f: \RR \to \RR$ such that $f'$ is surjective and nowhere flat
and for all $x_0 \in \RR$, if $f'(x_0) = -n^2$, $n \in \ZZ$,
then $(f''(x_0), f'''(x_0)) \ne (0,0)$.
A \textit{good nonlinearity} is an admissible nonlinearity
for which $f'(x_0) = -n^2$, $n \in \ZZ$, implies $f''(x_0) \ne 0$.
Most polynomials of even degree are good nonlinearities.

Let $\RR^2_{xy} \subset \RR^3$ be the plane $z = 0$ and
$I_n = [2\pi n - {\pi}/{2}, 2\pi n + {\pi}/{2}]$
for $n \in \ZZ$, $n > 0$.
Let $\cone_n$ be the cone
$x^2 + y^2 = \tan^2 z$, $z \in I_n$ and
$\Sigma = \RR^2_{xy} \cup \bigcup_{n > 0} \cone_n$.
The real separable infinite dimensional Hilbert space
will be denoted by $\HH$.


It turns out that the critical set $\Cc$ may include (at most countably many)
isolated points:
let $\Cc^\ast \subseteq \Cc$ be obtained from $\Cc$ by removing such points.

\begin{theo}
\label{theo:main}
Let $f: \RR \to \RR$ be an admissible nonlinearity:
the pair $(H^p(\Ss^1), \Cc^\ast)$ is diffeomorphic to the pair
$(\RR^3, \Sigma) \times \HH$;
if $f$ is a good nonlinearity then $\Cc^\ast = \Cc$.
\end{theo}

The present paper can also be considered a continuation of \cite{BST2},
where Theorem \ref{theo:main} is proved for the \textit{linear case}
$f(x) = x^2/2$:
the phrasing is justified by the fact that $f'(x) = x$
and therefore $DF(u)v = - v'' + uv$.
In this case we trivially have $\Cc^\ast = \Cc$. 
Notice that in the
periodic case, unlike the Dirichlet case, the critical set $\Cc$ has
singular points and is not a Hilbert manifold. A model of $\Cc$
at singular points is obtained by the study of the {\it monodromy map}.
We now review the notation and results used in
\cite{BST2} for the linear case, which will be heavily used in this paper.

Let $\Pi: G = \widetilde{SL(2,\RR)} \to SL(2,\RR)$
be the universal cover of $SL(2,\RR)$.
The elements of $SL(2,\RR)$ and $G$ will be referred to
as matrices and lifted matrices, respectively.
Recall that $G$ is itself a Lie group diffeomorphic to $\RR^3$.
For a potential $q \in H^p(\Ss^1)$,
let $v_1, v_2 \in H^{p+2}([0,2\pi])$ be the fundamental solutions
\[ v_i''(t) = q(t) v_i(t),\quad
v_1(0) = 1, \; v_1'(0) = 0, \; v_2(0) = 0, \; v_2'(0) = 1 \]
and define the {\it lifted fundamental matrix}
$\tilde\Phi: [0,2\pi] \to G$ by $\tilde\Phi(0) = I$ and
\[ \Pi(\tilde\Phi(t)) = \begin{pmatrix}
v_1(t) & v_2(t) \\ v_1'(t) & v_2'(t) \end{pmatrix}. \]
The \textit{monodromy map} $\mu: H^p(\Ss^1) \to G$
is the lifting $\mu(q) = \tilde\Phi(2\pi)$;
the projection $\Pi(\mu(q)) \in SL(2,\RR)$
is the so called Floquet multiplier associated to the potential $q$.
Notice that $L^1([0,2\pi])$ is another legitimate domain for $\mu$.

It is easy to verify that $u \in \Cc$
if and only if $\mu(u) \in T_2 \subset G$
where $T_2$ is the set of lifted matrices
of trace equal to $2$ (or, equivalently, with double eigenvalue $1$);
here $\tr(g) = \tr(\Pi(g))$ for a lifted matrix $g \in G$.
The image of $\mu$ is an open set $G_0 \subset G$ diffeomorphic to $\RR^3$.
Theorem 3 in \cite{BST2} constructs an explicit smooth diffeomorphism
$\Psi: G_0 \times \HH \to H^p(\Ss^1)$
such that $\mu \circ \Psi$ is the projection on the first coordinate.
It follows that $(H^p(\Ss^1), \Cc)$ is diffeomorphic via $\Psi^{-1}$
to $(G_0, T_2 \cap G_0) \times \HH$.

In the general case,
we define the \textit{nonlinear monodromy map} $\mu_f: H^p(\Ss^1) \to G_0$ by
$\mu_f(u) = \mu(f' \circ u)$.
Unlike $\mu$, the map $\mu_f$ cannot in general
be extended to $L^1([0,2\pi])$ but
$L^\infty([0,2\pi])$ will be enough for the purposes of this paper.
For admissible nonlinearities $f$,
it turns out that the image of $\mu_f$ equals $G_0$ and
we still have that $\Cc = \mu_f^{-1}(T_2 \cap G_0)$.
We shall not construct a counterpart of the diffeomorphism $\Psi$
for the general nonlinearity:
instead, we prove the contractibility of fibers of $\mu_f$
and show that this information suffices to model $\Cc$.

The local behavior of $\mu_f$
is particularly nasty at constant functions $u$;
on the other hand, such functions form
a subspace of infinite codimension and can therefore be excised
without changing the homotopy type of the domain.
Set $X^{\ast} = H^p(\Ss^1) \smallsetminus \{ u \textrm{ constant} \}$:
as we shall see in Theorem \ref{theo:contract},
the map $\mu_f: X^{\ast} \to G_0$ is a surjective submersion
with contractible fibers.
In order to complete the proof of Theorem \ref{theo:main}
we need some results in infinite dimensional topology.
More precisely, we present in Theorem \ref{theo:danA}
a normal form near $\Cc$,
a Hilbert submanifold with singularities:
such singularities arise from the fact that $T_2 \subset G$
is diffeomorphic to a countable union of cones.

As in \cite{BST1}, \cite{MST2}, \cite{MST1},
a key ingredient is the contractibility of
level sets of certain functionals defined on infinite dimensional
spaces. The reader may see little in common among the several proofs.
A unifying feature is
that we first construct a fake homotopy and then fix it: it helps
that the functional can actually be extended to a larger infinite
dimensional space with a weaker topology. Theorem 2 in \cite{BST1},
transcribed below, allows for moving from one space
to another without changing the homotopy type of level sets.

\begin{theo}
\label{theo:th2bst}
Let $X$ and $Y$ be separable Banach spaces.
Suppose $i: Y \to X$ is a bounded, injective linear map
with dense image and $M \subset X$ is a smooth, closed submanifold
of finite codimension. Then $N = i^{-1}(M)$ is a smooth
closed submanifold of $Y$ and the restrictions
$i: Y - N \to X - M$ and $i: (Y,N) \to (X,M)$
are homotopy equivalences.
\end{theo}

Section 2 contains basic facts about the linear monodromy map $\mu$,
including a study of the effect of adding bumps $\ell_i$ to a potential $q$
as controlled perturbations of $\mu(q + \sum a_i \ell_i)$.
In Section 3 we compute the derivative of the nonlinear monodromy
map $\mu_f$ and extend the results for bumps to this case.
We also verify that under rather general hypothesis on $f$
the image of $\mu_f$ is $G_0$,
as in the linear case (Proposition \ref{prop:mufsurjective}).
The argument runs as follows:
for $g \in G_0$, we first obtain a discontinuous $u_0 \in L^\infty$
which is smoothened out yielding $u_1$ with $\mu_f(u_1) \approx g$;
the error is then corrected by adding appropriate bumps to $u_1$.
In Section 4 we prove that level sets of $\mu_f$ in $X^{\ast}$
are contractible (Theorem \ref{theo:contract}) by constructing homotopies:
we first obtain a fake homotopy by composing
a homotopy for the linear case with a (discontinuous!) right inverse for $f'$.
Smoothening and correction are then similar to that in the proof of
Proposition \ref{prop:mufsurjective}.
Section 5 contains the necessary results in infinite dimensional topology
(including Theorem \ref{theo:danA})
which take into account the presence of cone-like objects.
The ingredients are combined in Section 6
to complete the proof of Theorem \ref{theo:main};
we express our thanks to S. Ferry for providing first the arguments and
afterwards directing us to the appropriate reference \cite{Michael} for
what we call here Michael's theorem (Fact \ref{fact:michael}).

The second and third authors received the support of
CNPq, CAPES and FAPERJ (Brazil).
The second author acknowledges the kind hospitality of
The Mathematics Department of The Ohio State University
during the winter quarter of 2004.
We thank the referee for a very careful reading.

\vfil

%% file: periodicA1.tex
\section{The monodromy map $\mu$}

We begin with a recollection of elementary facts about
the second order ODE $v'' = q(t) v$ and
the universal cover $\Pi: G = \widetilde{SL(2,\RR)} \to SL(2,\RR)$.
For details,
including a careful tracking of differentiability classes, see \cite{BST2}.

\begin{figure}[ht]
\begin{center}
\psfrag{dG0}{$\partial G_0$}
\psfrag{T2}{$T_2$}
\epsfig{height=44mm,file=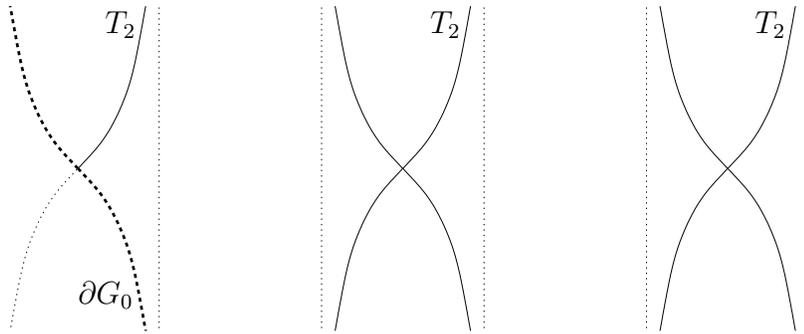}
\end{center}
\caption{The pair $(G_0, T_2 \cap G_0)$.}
\label{fig:g0p}
\end{figure}

The {\it left Iwasawa decomposition} is
the diffeomorphism $\phi_L: \RR \times (0,\infty) \times \RR \to G$
given by $\phi_L(0,1,0) = I$ and
\[
(\Pi \circ \phi_L)(\theta,\rho,\nu) =
\begin{pmatrix} \cos\theta & \sin\theta \\
-\sin\theta & \cos\theta \end{pmatrix}
\begin{pmatrix} \sqrt{\rho} & 0 \\ 0 & 1/\sqrt{\rho} \end{pmatrix} 
\begin{pmatrix} 1 & 0 \\ \nu/2 & 1 \end{pmatrix}.
\]
The set $G_0$ is defined by
$G_{0} = \phi_L((0,+\infty) \times (0,+\infty) \times \RR) \subset G$.
Equivalently, $g \in G_{0}$ if and only if
the variation in argument from $e_2$ to $g e_2$ is negative
(the variation in argument is computed along a path
$\gamma: [0,1] \to G$ joining $\gamma(0) = I$ to $\gamma(1) = g$).
The pair $(G_0, T_2 \cap G_0)$
is shown schematically in Figure \ref{fig:g0p}.
Here $T_2 \subset G$ is the set of lifted matrices of trace equal to $2$.
The thick dashed curve represents $\partial G_0$,
which is diffeomorphic to a plane;
the X's formed by crossing curves stand for the connected components of $T_2$,
diffeomorphic to cones with horizontal axis.
The dotted vertical lines represent the set of lifted matrices
of trace $0$ so that vertical regions contain matrices with trace
of alternating sign.
The pair $(\RR^3, \Sigma)$ constructed in the introduction
is diffeomorphic to $(G_0, T_2 \cap G_0)$:
$\Sigma = \RR^2_{xy} \cup \bigcup_{n > 0} \cone_n$
where $\RR^2_{xy} \subset \RR^3$ is the plane $z = 0$
and $\cone_n$ is the cone $x^2 + y^2 = \tan^2 z$, $|z - 2\pi n| < {\pi}/{2}$.

Given a potential $q \in L^\infty([0,2\pi])$,
let $v_1, v_2$ be the fundamental solutions
\[ v_i''(t) = q(t) v_i(t),\quad
v_1(0) = 1, \; v_1'(0) = 0, \; v_2(0) = 0, \; v_2'(0) = 1 \]
and define the {\it lifted fundamental matrix}
$\tilde\Phi: [0,2\pi] \to G_0$
and its projection $\Phi = \Pi \circ \tilde\Phi$
by $\tilde\Phi(0) = I$ and
\[ \Phi(t) = \Pi(\tilde\Phi(t)) = \begin{pmatrix}
v_1(t) & v_2(t) \\ v_1'(t) & v_2'(t) \end{pmatrix}. \]
Let $\vv: [0,2\pi] \to \RR^2$ be $\vv(t) = (v_1(t), v_2(t))$ 
so that
\[\vv''(t) = q(t) \vv(t),
\quad \vv(0) = e_1 = (1,0), \quad \vv'(0) = e_2 = (0,1). \]
The curve $\vv$ never passes through the origin so
its argument $\theta$ is well defined:
it is the only continuous function which satisfies
\[ \theta(0) = 0, \quad
\frac{\vv(t)}{||\vv(t)||} = (\cos \theta(t), \sin \theta(t)).\]
Since $\vv(t) \wedge \vv'(t) = 1$,
the function $\theta: [0, 2\pi] \to [0, \theta_M]$ is an increasing bijection.
Following \cite{BST2}, the \textit{orbit} $\rho: [0,\theta_M] \to (0,+\infty)$
associated to $\vv$ is
\[\rho(\theta(t)) = v_1^2(t) + v_2^2(t) = ||\vv(t)||^2. \]
For $\nu(\theta) = \rho'(\theta)/\rho(\theta)$,
it turns out that 
\( \phi_L( \theta(t), \rho(\theta(t)), \nu(\theta(t)) ) = \tilde\Phi(t) \).

Conversely, given an orbit, i.e.,
a number $\theta_M > 0$ and a function $\rho: [0,\theta_M] \to (0,+\infty)$
with $\rho(0) = 1$, $\rho'(0) = 0$ and
$\int_0^{\theta_M} \rho(\theta) d\theta = 2\pi$, consider
the parametrized curve $\sqrt{\rho(\theta)} (\cos\theta, \sin\theta)$.
There is a unique orientation preserving reparametrization
$\theta: [0,2\pi] \to [0,\theta_M]$ of the curve
so that equal areas around the origin are swept in equal times.
In other words, the reparametrized curve
\[ \vv(t) = \sqrt{\rho(\theta(t))} (\cos\theta(t), \sin\theta(t)) \]
satisfies $\vv(t) \wedge \vv'(t) = 1$ for all $t$.
Taking derivatives, $\vv(t) \wedge \vv''(t) = 0$ whence
$\vv''(t) = q(t) \vv(t)$ for some potential $q: [0,2\pi] \to \RR$.
The correspondence between potentials and orbits is called
the \textit{Kepler transform} in \cite{BST2}.

We now focus on the monodromy map $\mu: L^\infty([0,2\pi]) \to G_0$
given by $\mu(q) = \tilde\Phi(2\pi)$.
Recall that $\fg = sl(2,\mathbb R)$ is the space of $2\times 2$
real matrices $A = (a_{ij})$ with $\tr A=0$.
If $g\in G$ and $M \in T_e(G) = \fg$
we write $gM$ for the element in $T_g(G)$ obtained as the image of $M$ by the
differential of the translation $g$.
Alternatively, the natural identification $T_g(G) = T_{\Pi(g)}(SL(2,\RR))$
allows us to interpret $gM$ as the matrix product $\Pi(g) M$.
For an angle $\omega\in \mathbb R$, define $N_\omega \in \fg$ by 
\[ N_\omega =
\begin{pmatrix} -\sin\omega \cos\omega& -\sin^2\omega \\
\cos^2\omega & \sin\omega \cos\omega \end{pmatrix} =
\frac{1}{2} \begin{pmatrix} -\sin 2\omega & -1+\cos 2\omega \\
1+\cos 2\omega & \sin 2\omega \end{pmatrix}; \]
notice that the matrices $N_\omega$ form a circle
in the plane $a_{21} - a_{12}=1$.

\begin{prop}
\label{prop:Dxi}
The monodromy map $\mu: L^\infty([0,2\pi]) \to G_0$
is smooth with derivative given by
\[ (D\mu(q)) w =
\mu(q) \left( \int_0^{2\pi} w(t) \rho(\theta(t)) N_{\theta(t)} dt \right). \]
Furthermore, $\mu$ is uniformly continuous
on bounded subsets of $L^\infty([0,2\pi])$. 
\end{prop}

{\nobf Proof: }
The variation of the fundamental solutions $v_i$
with respect to the potential is a familiar computation (see \cite{PT}),
as are expressions for higher derivatives.
The formula above for the derivative of the monodromy is then easy.

An alternative way of obtaining this formula is to (temporarily) allow for
potentials which are distributions.
A straightforward computation yields
\[ \rho(\theta(t)) N_{\theta(t)} = (\Phi(t))^{-1} N_0 \Phi(t) =
\begin{pmatrix} -v_1(t) v_2(t) & -v_2^2(t) \\
v_1^2(t) & v_1(t) v_2(t) \end{pmatrix}. \]
Thus, if $\delta_{t_0}$ is the Dirac delta
centered at $t_0$, $0 < t_0 < 2\pi$, it is easy to verify that
\[ \mu(q + a\delta_{t_0}) =
\mu(q) (\tilde\Phi(t_0))^{-1} \begin{pmatrix}1 & 0 \\ a & 1\end{pmatrix}
\tilde\Phi(t_0) \]
and therefore $(D\mu(q)) \delta_{t_0} = \mu(q) \rho(\theta(t)) N_{\theta(t)}$.

Clearly, a potential $q \in L^\infty$ induces fundamental solutions,
monodromy matrices and functions $\theta$ and $\rho$ which are bounded
by simple expressions in $|q|_{L^\infty}$. The formula obtains
bounds for the derivative and therefore uniform continuity
on $L^\infty$-bounded sets.
\qed

Assume $p \ge 1$ so that $H^p([0,2\pi]) \subset C^0([0,2\pi])$.
The formula for $D\mu$ in the proof above implies the known fact
that the functions $v_1^2$, $v_1v_2$ and $v_2^2$ are taken
by $D\mu$ to a basis of $T_{\mu(q)}G$, whence $\mu$ is a submersion.
We need a more workable triple of generators, however.
A {\it bump} centered at $t_0 \in \Ss^1$ is a smooth nonnegative
function from $\Ss^1$ to $\RR$ whose support is a small interval centered
at $t_0$; the size of the support is the {\it width} of the bump.
Potentials will be altered by adding bumps
in order to adjust the value of $\mu$.


\begin{lemma}
\label{lemma:estimate}
Let $q \in H^p([0,2\pi])$ be a potential.
\begin{enumerate}[(a)]
\item{ Set $\epsilon^{-2} = |q|_{L^{\infty}}$.
If $0 < t^+ - t^- < \epsilon$ then
$\theta(t^-) < \theta(t^+) < \theta(t^-) + \pi$.}
\item{ Let $\ell_i$, $i = 1, 2, 3$, be bumps with disjoint supports
contained in an interval $[t^-, t^+]$ with
$\theta(t^+) < \theta(t^-) + \pi$.
Then the vectors $D\mu(q) \ell_i$, $i = 1, 2, 3$, are linearly independent.
In particular, $\mu: H^p([0,2\pi]) \to G_0$ and $\mu: H^p(\Ss^1) \to G_0$
are submersions.}
\end{enumerate}
\end{lemma}

{\nobf Proof: }
With the hypothesis of item (a), let $v$ be the solution of
$v(t^-) = 0$, $v'(t^-) = 1$, $v''(t) = q(t) v(t)$;
we claim that $v'(t) > 0$ for all $t \in (t^-,t^+)$.
Assume by contradiction that $t_{\max} \le t^+$ is the first
local maximum of $v$ (after $t^-$) so that $v'(t_{\max}) = 0$;
set $v_{\max} = v(t_{\max})$.
The maximum value of $v'$ in the interval $[t^-,t_{\max}]$ is
at least $v_{\max}/\epsilon$ and therefore there exists $t$
in this interval with $v''(t) < -v_{\max}/(\epsilon^2)$ and
therefore $v''(t) = q(t) v(t) < -\epsilon^{-2} v(t)$ and
$|q(t)|_{L^{\infty}} > \epsilon^{-2}$, a contradiction.

The claim implies the linear independence of
the vectors $(v(t^-), v'(t^-))$ and $(v(t^+), v'(t^+))$.
Since $v$ is in the linear span of $v_1$ and $v_2$,
$(v_1(t^-), v_2(t^-))$ and $(v_1(t^+), v_2(t^+))$
are not collinear, concluding the proof of item (a).

\begin{figure}[ht]
\begin{center}
\vskip6mm
\psfrag{th1m}{$\theta(t_1^-)$}
\psfrag{th1p}{$\theta(t_1^+)$}
\psfrag{th2m}{$\theta(t_2^-)$}
\psfrag{th2p}{$\theta(t_2^+)$}
\psfrag{th3m}{$\theta(t_3^-)$}
\psfrag{th3p}{$\theta(t_3^+)$}
\epsfig{height=50mm,file=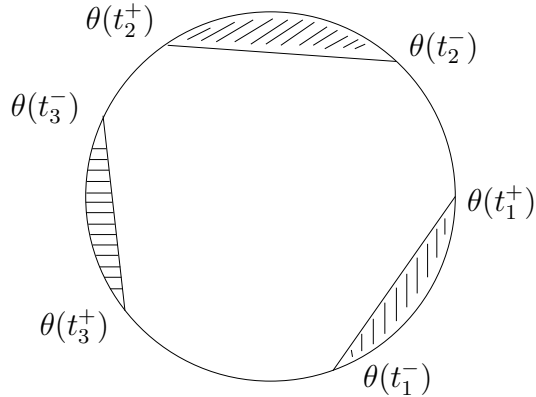}
\end{center}
\caption{Convex hulls of three arcs}
\label{fig:3arcs}
\end{figure}

Let $[t_i^-,t_i^+]$ be the support of $\ell_i$.
Without loss of generality,
\[ \theta(t^-) \le \theta(t_1^-) < \theta(t_1^+) <
\theta(t_2^-) < \theta(t_2^+) <
\theta(t_3^-) < \theta(t_3^+) \le \theta(t^+) < \theta(t^-) + \pi. \]
Note that $\theta (t)$ is strictly increasing in $t$.
Write
\[ (\mu(q))^{-1} (D\mu(q)) \ell_i =
\int_{t_i^-}^{t_i^+}
\ell_i(t) \rho(\theta(t)) N_{\theta(t)} dt \notag
=
\int_{\theta(t_i^-)}^{\theta(t_i^+)}
\frac{\ell_i(\theta^{-1}(\tau)) \rho(\tau)}{\theta'(\theta^{-1}(\tau))}
N_{\tau} d\tau. \]
Thus, up to a positive multiplicative factor,
$(\mu(q))^{-1} (D\mu(q)) \ell_i$ is a convex combination of matrices
$N_\tau$ in the arc from $N_{\theta(t_i^-)}$ to $N_{\theta(t_i^+)}$
and, in particular, lies in the plane $a_{21} - a_{12} = 1$.
Figure \ref{fig:3arcs} illustrates that if we take a point in the convex
hull of each arc we necessarily form a non-degenerate triangle
on this plane and therefore the vectors $(\mu(q))^{-1} (D\mu(q)) \ell_i$
are linearly independent.
\qed

\bigskip
\section{The nonlinear monodromy map $\mu_f$}

For a smooth nonlinearity $f \in C^\infty(\RR,\RR)$
define the nonlinear monodromy map $\mu_f: L^\infty([0,2\pi]) \to G_0$ by
$\mu_f(u) = \mu(f' \circ u)$.
We often consider restrictions of $\mu_f$ to smaller spaces
such as $\mu_f|_{H^p(\Ss^1)}$, $p \ge 1$.
The critical set $\Cc$ of the operator $F: H^p(\Ss^1) \to H^{p-2}(\Ss^1)$
defined by $F(u) = - u'' + f(u)$ is
$\Cc = (\mu_f|_{H^p(\Ss^1)})^{-1}(T_2 \cap G_0)$,
where $T_2 \subset G$ is the set of lifted matrices with trace $2$.
A real separable Banach space $X$ is \textit{smoothing}
if the inclusions
$C^\infty(\Ss^1) \subset X \subset L^\infty([0,2\pi])$ are continuous.

\goodbreak

\begin{prop}
\label{prop:localmu}
Let $f$ be a smooth nonlinearity.
\begin{enumerate}[(a)]
\item{ The map $\mu_f: L^\infty([0,2\pi]) \to G_0$ is smooth.
Given $u \in L^\infty(\Ss^1)$,
define $\theta$ and $\rho$ as before for the potential $f' \circ u$.
Then, for any $v \in L^\infty([0,2\pi])$,
\[ (D\mu_f(u)) v = \mu_f(u) \left(
\int_0^{2\pi} f''(u(t)) v(t) \rho(\theta(t)) N_{\theta(t)} dt \right). \]}
\item{ Let $X$ be a smoothing Banach space; assume $f'$ nowhere flat.
If $u \in X$ is not constant or is constant with $f''(u) \ne 0$,
then the derivative $D(\mu_f|_{X})(u)$ is a surjective linear map and
$\mu_f|_{X}$ is a local submersion at $u$.
If $u \in X$ is constant which is a local extremum of $f'$ then
$u$ is an isolated point in the level set
$(\mu_f|_{X})^{-1}(\mu_f(u))$.}
\end{enumerate}
\end{prop}

{\nobf Proof: }
Item (a) follows directly from Proposition \ref{prop:Dxi}.
The first case in item (b) follows from item (a).
From standard oscillation theory, if $f'(u_1(t)) \le f'(u_2(t))$,
$u_1 \ne u_2$,
then the arguments of either column of $\mu_f(u_1)$ are respectively larger
than those of $\mu_f(u_2)$.
This implies the rest of (b), which alternatively
follows from a simple local argument using Proposition \ref{prop:Dxi},
left to the reader.
\qed

If $f''$ and $f'''$ have no common zeroes then these
cover all the possibilities for $u \in X$.
For flat $f'$, more complicated scenarios might occur,
including accumulations of isolated points in level sets
$(\mu_f|_{X})^{-1}(g)$.
Constant functions $u$ are therefore potentially nasty objects:
for a smoothing Banach space $X$,
let $X^{\ast} = X - \{u = \textrm{const.}\}$:
recall that $X^{\ast}$ is homeomorphic to the separable Hilbert space $\HH$
and diffeomorphic to $\HH$ if $X$ itself is Hilbert.
Denote by $\lambda$ the Lebesgue measure in $\RR$.

\begin{prop}
\label{prop:submersion}
Let $X$ be a smoothing Banach space and
$f: \RR \to \RR$ be a smooth function with nowhere flat derivative.
Let $u \in X$ and $\ell_i$, $i = 1, 2, 3$, be bump functions
with disjoint supports contained in an interval $I \subset (0,2\pi)$
such that $f''(u(t)) \ne 0$ for all $t \in I$.
Then, for every $E > 0$ there exists $\epsilon > 0$ such that
if $|u|_{L^\infty} < E$ and $\lambda(I) < \epsilon$
then the vectors $D\mu_f(u) \ell_i$, $i = 1, 2, 3$, are linearly independent.
Thus, the function $\mu_f: X^{\ast} \to G_0$ is a submersion.
\end{prop}

Notice that we are not claiming (yet) that $\mu_f: X^{\ast} \to G_0$
is surjective: this is Proposition \ref{prop:mufsurjective}.

{\noindent\bf Proof: }
Let $E'$ be such that $|x| < E$ implies $|f'(x)| < E'$
and take $\epsilon = \sqrt{1/E'}$.
We may assume without loss of generality
that $f''(u(t)) > 0$ for all $t \in I$ so that
$\tilde\ell_i(t) = f''(u(t)) \ell_i(t)$, $i = 1, 2, 3$,
are bump functions.
From lemma \ref{lemma:estimate},
the three vectors $D\mu(f' \circ u)\tilde\ell_i = D\mu_f(u) \ell_i$
are linearly independent and we are done.
\qed

The following corollary generalizes the first part of
Proposition \ref{prop:submersion};
it will be needed in the proof of Theorem \ref{theo:contract}.

\begin{coro}
\label{coro:submersionX}
Let $X \subset C^0([0,2\pi])$ be a separable Banach space with continuous
inclusion such that the closure of $X$ is either
$C^0(\Ss^1)$ or $C^0([0,2\pi])$.
Let $f: \RR \to \RR$ be a smooth function with nowhere flat derivative.
The function $\mu_f: X^{\ast} \to G_0$ is a submersion.
\end{coro}

{\noindent\bf Proof: }
We must prove that $D\mu_f(u): X \to T_{\mu_f(u)} G$ is surjective.
By Proposition \ref{prop:submersion},
$D\mu_f(u)(C^0([0,2\pi])) = D\mu_f(u)(C^0(\Ss^1)) =  T_{\mu_f(u)} G$
and the result follows by density.
\qed

The rest of the section is dedicated to showing in
Proposition \ref{prop:mufsurjective}
that $\mu_f: X \to G_0$ is surjective for smoothing spaces $X$.
Notice that the image of $\mu_f$ clearly equals $G_0$ if
$f': \RR \to \RR$ is a diffeomorphism.
Indeed, given $g \in G_0$
let $q \in C^\infty(\Ss^1)$ with $\mu(q) = g$:
the function $u = (f')^{-1} \circ q$ satisfies
$\mu_f(u) = \mu(f' \circ u) = g$.
Similarly, Theorem \ref{theo:main} follows easily from 
results in \cite{BST2} if $f'$ is a diffeomorphism.
The idea in the proofs of Proposition \ref{prop:mufsurjective} and
Theorem \ref{theo:contract} is to first operate with a discontinuous
right inverse for $f'$, then fix the discontinuities of the potential
and add bumps in order not to change the monodromy.

We call a function $h: \RR \to \RR$ {\it piecewise smooth} if
there exists a discrete set $Y_h$ such that
if $I$ is a connected component of $\RR - Y_h$
then there exists a continuous function $h_I: \overline I \to \RR$
such that $(h_I)|_I = h|_I$ is smooth.
Notice that an element $x \in Y_h$ may be a discontinuity of $h$
or a point where $h$ is continuous but not smooth,
such as $x = 0$ for $h(x) = x^{1/3}$.
The proof of the following lemma is left to the reader.

\begin{lemma}
\label{lemma:inverse}
If $h: \RR \to \RR$ is smooth, nowhere flat and surjective
then there exists a piecewise smooth function $h^\sharp: \RR \to \RR$
with $(h \circ h^\sharp)(x) = x$ for all $x$.
\end{lemma}

Discontinuities in $u_0 = (f')^\sharp \circ q$ will be handled
by considering a smooth function $u_1$ such that $|u_1 - u_0|_{L^1}$ is small.
In particular, $\mu_f(u_1)$ is close to $\mu_f(u_0)$:
adding appropriate bumps to $u_1$ produces a smooth function $u_2$
with $\mu_f(u_2) = \mu_f(u_0) = g$.

Endow the Lie algebra $\fg = sl(2,\RR)$ of $G$ with an Euclidean metric
and let $\expg: \fg \to G$ be the exponential map:
take $\epsilon_\fg > 0$ to be such that the restriction of $\expg$
to the ball in $\fg$ of center $0$ and radius $\epsilon_\fg$
is injective.
For $\epsilon \le \epsilon_\fg$, let $B_\epsilon \subset G$
be the image under $\expg$ of the ball of center $0$ and radius $\epsilon$.
The sets $B_\epsilon$ are invariant under inversion
($h \in B_\epsilon$ if and only if $h^{-1} \in B_\epsilon$)
but not under conjugation:
in general, $B_\epsilon \ne h^{-1} B_\epsilon h$.
By continuity, given $\epsilon > 0$, $\epsilon \le \epsilon_\fg$,
there exists $\delta > 0$ such that $B_\delta B_\delta \subseteq B_\epsilon$
(i.e., if $h_1, h_2 \in B_\delta$ then $h_1h_2 \in B_\epsilon$).

\begin{lemma}
\label{lemma:LinftyBeps}
Consider the monodromy $\mu_f: L^\infty([0,2\pi]) \to G_0$
for a smooth nonlinearity $f: \RR \to \RR$.
Then for all $M \in \RR$ and for all $\epsilon > 0$
there exists $\delta > 0$ such that
for all $u_1, u_2 \in L^\infty([0,2\pi])$
\[ |u_1|_{L^\infty}, |u_2|_{L^\infty} < M, \quad
\lambda(\{t| u_1(t) \ne u_2(t)\}) < \delta \quad
\Longrightarrow \quad
\mu_f(u_2) (\mu_f(u_1))^{-1} \in B_\epsilon. \]
\end{lemma}

{\nobf Proof: }
Let $\tilde M = \sup(\{|f'(x)|, |x| < M\})$.
Then $|f'\circ u_i|_{L^\infty} < \tilde M$
and the condition $\lambda(\{t| u_1(t) \ne u_2(t)\}) < \delta$ implies
$|f'\circ u_2 - f'\circ u_1|_{L^1} < 2\tilde M\delta$.
From Proposition \ref{prop:Dxi},
the monodromy map $\mu: L^1([0,2\pi]) \to G_0$
is uniformly continuous with respect to the $L^1$ norm
in the set $\{ q \in L^1([0,2\pi]), |q|_{L^\infty} < \tilde M \}$,
completing the proof.
\qed

It will be convenient to restrict the monodromy to intervals.
For $T = [t_-, t_+] \subset [0, 2\pi]$ and $u \in L^\infty(T)$,
let $\tilde\Phi: T \to G$ be the only solution of
\[ \tilde\Phi(t_-) = I,\quad \tilde\Phi'(t) =
\begin{pmatrix} 0 & 1 \\ f'(u(t)) & 0 \end{pmatrix} \tilde\Phi(t) \]
and define the {\it $T$-variation} $\mu_{f,T}(u) = \tilde\Phi(t_+)$.
Variations juxtapose, and the order is important:
if $t_0 < t_1 < t_2$ and $u \in L^\infty([t_0,t_2])$
then $\mu_{f,[t_0,t_2]}(u) =
\mu_{f,[t_1,t_2]}(u) \mu_{f,[t_0,t_1]}(u)$.
In particular,
$\mu_f(u) = \mu_{f,T_{2}}(u) \mu_{f,T_{1}}(u)$.

\begin{prop}
\label{prop:mufsurjective}
Let $X$ be a smoothing Banach space
and $f: \RR \to \RR$ be a smooth function such that
$f'$ is nowhere flat and surjective.
Then $\mu_f: X^\ast \to G_0$ is surjective.
\end{prop}

{\nobf Proof: }
Take $g \in G_0$.
For $y_0$ a regular point of $(f')^\sharp$, set $x_0 = (f')^\sharp(y_0)$
so that $x_0$ is a regular point of $f'$ and $f'(x_0) = y_0$.
From the Kepler transform construction, outlined in the previous section,
there exists $\epsilon_1 > 0$ and $q \in C^\infty([0,2\pi])$
such that $\mu(q) = g$, $q$ is nonflat in $(\epsilon_1, 2\pi - \epsilon_1)$
and $q$ is constant equal to $y_0$ in $[0,\epsilon_1] \cup [\epsilon_1,2\pi]$.
Take $u_0 = (f')^\sharp \circ q$:
clearly $u_0 \in L^\infty([0,2\pi])$, $\mu_f(u_0) = g$.
The function $u_0$ is piecewise smooth:
let $T_{1} = [0,2\pi-\epsilon_1]$, $T_{2} = [2\pi - \epsilon_1, 2\pi]$
and $Y \subset T_{1}$ be the (discrete) set of discontinuities of $u_0$.
We need to alter $u_0$ in an open neighborhood $\breve Y \supset Y$,
$\breve Y \subset T_1$,
and add bumps in $T_2$ so as to obtain $u_2 \in X$ with $\mu_f(u_2) = g$.

Let $\ell_i$, $i = 1, 2, 3$, be bumps with disjoint supports contained
in $T_{2}$. From Lemma \ref{lemma:estimate} we may assume
the linear independence of the vectors $D\mu(q)\ell_i$
and therefore of the vectors $D\mu_{f,T_{2}}(u_0)\ell_i$.
From the inverse function theorem,
there exists an open neighborhood $B \subset G_0$ of $\mu_{f,T_{1}}(u_0)$
such that for any $h \in B$ there exist $a_i$, $i = 1,2,3$,
which adjust the monodromy:
\[ \mu_{f,T_{2}}\left( u_0+ \sum_{i = 1,2,3} a_i \ell_i \right) h
= g. \]
Let $E = 2 |u_0|_{L^\infty}$.
Use Lemma \ref{lemma:LinftyBeps} to obtain $\epsilon_2 > 0$ be such that
\[ |u_1|_{L^\infty} < 2E, \quad
\lambda(\{t \in T_{1} | u_0(t) \ne u_1(t)\}) < \epsilon_2
\quad\Rightarrow\quad
\mu_{f,T_{1}}(u_1) \in B. \]
Choose an open neighborhood $\breve Y \supset Y$ with
$\breve Y \subset T_1$ and $\lambda(\breve Y) < \epsilon_2$.
Let $u_1$ be an arbitrary smooth function coinciding with
$u_0$ in $[0,2\pi] \smallsetminus \breve Y$, $|u_1|_{L^\infty} < 2E$.
We have $h = \mu_{f,T_1}(u_1) \in B$ and therefore
there exist $a_i$, $i = 1,2,3$, such that 
$u_2 = u_1 + \sum a_i \ell_i$ satisfies $\mu_f(u_2) = g$.
\qed

\bigskip

\section{Levels of generic $\mu_f$ are contractible}

For the linear case, the levels sets of the monodromy map $\mu$
were explicitly parametrized in \cite{BST2}.
Indeed, for a smoothing Banach space $X$,
let $\Lev^X(g) \subset X$ be the level set $\mu^{-1}(\{g\})$,
where $\mu: X \to G_0 \subset G$ is the monodromy map.
Theorem $3$ in \cite{BST2} gives a diffeomorphism between
$\Lev^X(g)$ and $\HH$ for $X = H^p([0,2\pi])$, $p \ge 0$:
in particular, $\Lev^X(g)$ is contractible.
In this section we generalize this last result:
let $\Lev^X_f(g)$ be the level set
$\mu_f^{-1}(\{g\}) \cap X^\ast$ where
$X^{\ast} = X - \{u = \textrm{const.}\}$.

\begin{theo}
\label{theo:contract}
Let $X$ be a smoothing Banach (resp. Hilbert) space and
$f: \RR \to \RR$ be a smooth function such that
$f'$ is nowhere flat and surjective.
Then for any $g \in G_0$ the level set $\Lev^X_f(g) \subset X^{\ast}$ is
a nonempty contractible Banach (resp. Hilbert) submanifold of codimension $3$.
\end{theo}


The fact that $\Lev^X_f(g)$ is a submanifold
follows from Proposition \ref{prop:submersion}.
As discussed in Proposition \ref{prop:localmu},
the more natural level set $\mu_f^{-1}(\{g\}) \subset X$
may not be a Banach manifold at
constants $x$ for which $f''(x) = 0$.
Removing such points from $\mu_f^{-1}(\{g\})$ one obtains
a Banach manifold $\Lev_0$.
Notice that in the definition of $\Lev^X_f(g)$
we remove from $\mu_f^{-1}(\{g\})$ all constant functions
so that $\Lev^X_f(g) \subseteq \Lev_0 \subseteq \mu_f^{-1}(\{g\})$.
The inclusion $\Lev^X_f(g) \subseteq \Lev_0$
is a homotopy equivalence since the two Banach manifolds
differ by the excision of a subset of infinite codimension;
the set $\mu_f^{-1}(\{g\})$, on the other hand,
may contain additional isolated points.

It might seem that different choices of $X$ would demand different arguments.
Fortunately, this is not so: Theorem 2 from \cite{BST1},
transcribed in the introduction,
implies the homotopy equivalence of the spaces $\Lev^X_f(g)$
for different choices of $X$.
It is then sufficient to prove that all homotopy groups of
$\Lev^X_f(g)$ are trivial.
For the standard monodromy map $\mu$, i.e.,
for the linear case $f(x) = x^2/2$,
this result is proved in \cite{BST2}.
Now, for $X = H^1(\Ss^1)$, given a loop $\gamma: \Ss^k \to \Lev^{H^1}_f(g)$,
define $\gamma_\textrm{Lin}: \Ss^k \to H^1(\Ss^1)$
by $\gamma_\textrm{Lin}(s) = f' \circ \gamma(s)$ so that
$\mu(\gamma_\textrm{Lin}(s)) = g$ for all $s \in \Ss^k$.
From \cite{BST2}, the loop $\gamma_\textrm{Lin}$
admits an extension $\Gamma_\textrm{Lin}: \BB^{k+1} \to H^1(\Ss^1)$ 
with $\mu(\Gamma_\textrm{Lin}(s)) = g$ for all $s \in \BB^{k+1}$.
In order to obtain $\Gamma: \BB^{k+1} \to \Lev^{H^1}_f(g)$ 
such that $\mu_f(\Gamma(s)) = g$, we might want to define
$\Gamma_\textrm{Lin}(s) = f' \circ \Gamma(s)$:
such a construction will not respect $H^1(\Ss^1)$
if $f'$ is not invertible.
As in Proposition \ref{prop:mufsurjective},
the idea of the proof is to use the piecewise smooth right inverse
$(f')^\sharp$ of $f'$:
we must now regularize a family of discontinuous potentials
without changing their monodromies.

{\nobf Proof of Theorem \ref{theo:contract}: }
For $\epsilon > 0$, let $R_\epsilon \subset \CC$
be the closed rectangle of complex numbers $z = a+bi$,
$a \in [-\epsilon, 2\pi+\epsilon]$, $b \in [-\epsilon, \epsilon]$
and let $A_\epsilon$ be the set of continuous functions
$f: R_\epsilon \to \CC$ which are holomorphic in the interior of $R_\epsilon$
and which satisfy $f(t) \in \RR$ for $t \in [-\epsilon, 2\pi+\epsilon]$
and $f(z+2\pi) = f(z)$ whenever $z$ and $z+2\pi$ are in $R_\epsilon$.
Clearly, the inclusion $i: A_\epsilon \to H^p(\Ss^1)$ is bounded and,
from Theorem \ref{theo:th2bst} (with $X = H^p(\Ss^1)$ and $Y = A_\epsilon$)
and Corollary \ref{coro:submersionX},
we may assume $\gamma: \Ss^k \to \Lev^{A_\epsilon}_f(g)$.
From the same theorem, it suffices to extend $\gamma$ to
$\Gamma: \BB^{k+1} \to \Lev^{C^0(\Ss^1)}_f(g)$.

Since $\gamma(s)$ is analytic for all $s$, each $\gamma(s)$ is nowhere flat.
For a positive integer $r$, a smooth function $g: \RR \to \RR$
is {\it nowhere $r$-flat} if, for all $x \in \RR$,
there exists $r' = r'(x) \in \ZZ$, $0 < r' \le r$,
such that $g^{(r')}(x) \ne 0$.
Continuity in the $A_\epsilon$ norm yields a uniform bound:
there exists $r_\gamma$ such that each $\gamma(s)$ is nowhere $r_\gamma$-flat
(this is the reason we introduced the space $A_\epsilon$).
Assume $f'$ to be nowhere $r_f$-flat; set $r = r_\gamma r_f$.
Let $\gamma_\textrm{Lin}: \Ss^k \to \Lev^{C^{r}(\Ss^1)}(g)$ be defined
by $(\gamma_\textrm{Lin}(s))(t) = f'((\gamma(s))(t))$;
clearly, $\gamma_\textrm{Lin}(s)$ is nowhere $r$-flat for all $s \in \Ss^k$.
Again by Theorem \ref{theo:th2bst}, now setting $X = C^{r}(\Ss^1)$
and $Y = A_\epsilon$, there exists
$\gamma_\textrm{An}: \Ss^k \to \Lev^{A_\epsilon}(g)$
arbitrarily close to $\gamma_\textrm{Lin}$ and a homotopy
$\Gamma_{\textrm{Lin},\textrm{An}}: [0,1] \times \Ss^k \to \Lev^{C^{r}(\Ss^1)}$
joining $\gamma_\textrm{Lin}$ and $\gamma_\textrm{An}$.
We may furthermore assume that $\Gamma_{\textrm{Lin},\textrm{An}}(\tau,s)$
is nowhere ${r}$-flat for all $(\tau,s) \in [0,1] \times \Ss^k$.
Since $\Lev^{A_\epsilon}(g)$ is contractible there exists
$\Gamma_\textrm{An}: \BB^{k+1} \to \Lev^{A_\epsilon}(g)$
extending $\gamma_\textrm{An}$.
Juxtapose $\Gamma_{\textrm{Lin},\textrm{An}}$ and $\Gamma_\textrm{An}$
to define
$\Gamma_{\textrm{Lin}}: \BB^{k+1} \to \Lev^{C^{r}(\Ss^1)}(g)$
extending $\gamma_\textrm{Lin}$:
\[ \Gamma_\textrm{Lin}(s) = \begin{cases}
\Gamma_\textrm{An}(2s),&||s|| \le 1/2,\\
\Gamma_{\textrm{Lin},\textrm{An}}(2-2||s||,s/||s||),&||s|| \ge 1/2;
\end{cases} \]
notice that $\Gamma_{\textrm{Lin}}(s)$
is nowhere flat for all $s \in \BB^{k+1}$.
Define $\tilde\Gamma_\textrm{Lin}: \BB^{k+1} \times \Ss^1 \to \RR$
by $\tilde\Gamma_\textrm{Lin}(s,t) = (\Gamma_\textrm{Lin}(s))(t)$
(there will be similar correspondences between
$\Gamma_1$, $\Gamma_2$, $\Gamma_3$ below and their counterparts 
$\tilde\Gamma_1$, $\tilde\Gamma_2$, $\tilde\Gamma_3$;
deformations are often easier to describe from this second point of view).

\begin{figure}[ht]
\begin{center}
\epsfig{height=40mm,file=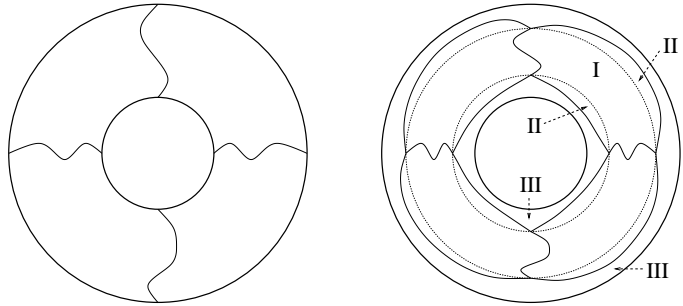}
\end{center}
\caption{The sets of discontinuities $Y_1$ and $Y_2$ of 
$\Gamma_1$ and $\Gamma_2$ (here, $k = 0$).}
\label{fig:Gamma2}
\end{figure}

Let $(f')^{\sharp}$ be as in lemma \ref{lemma:inverse}.
Let $\tilde\Gamma_1: \BB^{k+1} \times \Ss^1 \to \RR$
be the bounded (possibly discontinuous) function
$(f')^{\sharp} \circ \tilde\Gamma_\textrm{Lin}$:
define $\Gamma_1: \BB^{k+1} \to L^\infty([0,2\pi])$
by $(\Gamma_1(s))(t) = \tilde\Gamma_1(s,t)$.
Notice that $\mu_f(\Gamma_1(s)) = g$ for all $s \in \BB^{k+1}$ since
\[ \mu_f(\Gamma_1(s)) = \mu(f' \circ \Gamma_1(s)) =
\mu(\Gamma_\textrm{Lin}(s)) = g. \]
Let $Y_1 \subset \BB^{k+1} \times \Ss^1$ be the set of discontinuities
of $\tilde\Gamma_1$.
From nowhere-flatness, $Y_1$ intersects circles $(s,\cdot)$
in sets of measure zero; call this set $Y_1(s) \subset \Ss^1$:
\[ Y_1(s) = \Pi_{\Ss^1} (Y_1 \cap (\{s\} \times \Ss^1)). \]
Recall that there may exist points where $(f')^\sharp$
is continuous but not smooth: at such points, $\tilde\Gamma_1$
loses some differentiability but since our aim is to construct
$\Gamma: \BB^{k+1} \to \Lev^{C^0(\Ss^1)}_f(g)$,
this is not an issue and such points require no special treatment.
We modify $\Gamma_1$ by adding a thick shell to its domain
so that discontinuities will stay away from the boundary of the domain.
With the identification $\Ss^1 = \RR/(2\pi\ZZ)$, set $t_0  = 0 \in \Ss^1$;
define
\[ \tilde\Gamma_2(s,t) =
\begin{cases}
\tilde\Gamma_1(2s,t),&\textrm{if }||s|| \le 1/2,\\
\tilde\Gamma_1(s/||s||,t),&\textrm{if }||s|| \le
\frac{1}{2} + \frac{1}{10} d(t, (Y_1(s/||s||) \cup \{t_0\})),\\
\gamma(s/||s||,t),&\textrm{otherwise.}
\end{cases} \]
The three cases are indicated by I, II and III in the right side
of Figure \ref{fig:Gamma2}.
The function $\tilde\Gamma_2$ is also discontinuous:
let $Y_2 \subset \BB^{k+1} \times \Ss^1$ be its set of discontinuities.
As for $Y_1$, define $Y_2(s) \subset \Ss^1$ to be
the projection on $\Ss^1$ of the intersection of $Y_2$
with the circle $(s, \cdot)$;
clearly, the measure of $Y_2(s)$ is zero for all $s \in \BB^{k+1}$.
The closed set $Y_2$ keeps away from the boundary
(it is contained in a ball of radius $0.85$)
and $\tilde\Gamma_2$ is a (discontinuous) extension of $\tilde\gamma$.
Notice that $\mu_f(\Gamma_2(s)) = g$ for all $s$
since, for $s \le 1/2$, $\Gamma_2(s) = \Gamma_1(2s)$ and,
for $s \ge 1/2$,
\[ \mu_f(\Gamma_2(s)) = \mu(f' \circ \Gamma_2(s)) =
\mu(\gamma_\textrm{Lin}(s/||s||)) = g. \]
Assume $|\tilde\Gamma_2(s,t)| < E$ for all $(s,t) \in \BB^{k+1} \times \Ss^1$.
We must regularize $\tilde\Gamma_2$ without changing the monodromy:
at this point a sketch of what comes ahead is appropriate.

We will take an open set
$\breve Y_2 \subset \BB^{k+1} \times \Ss^1$, $Y_2 \subset \breve Y_2$
far from the boundary (i.e., $(s,t) \in \breve Y_2$ implies $||s|| < 9/10$)
and within it $\tilde\Gamma_2$ will be altered to obtain
a continuous function $\tilde\Gamma_3$
with $|\tilde\Gamma_3(s,t)| < 2E$ for all $(s,t) \in \BB^{k+1} \times \Ss^1$.
Notice that $\Gamma_3$ extends the original loop $\gamma$ continuously
but has slightly wrong values for the monodromy in the interior.

\goodbreak

\begin{figure}[ht]
\begin{center}
\epsfig{height=40mm,file=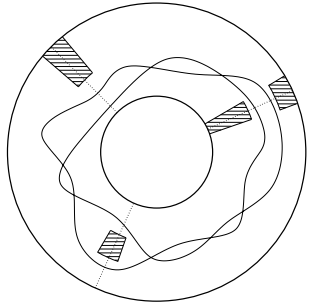}
\end{center}
\caption{The sets $K_j \times T_j \subset \BB^{k+1} \times \Ss^1$ avoid $Y_2$
(here, $k = 0$).}
\label{fig:boia}
\end{figure}

The monodromy of $\Gamma_3$ will be corrected with bumps
roughly as done in Proposition \ref{prop:mufsurjective}:
in this parametrized version, however, we have to find
places to support the bumps (the shaded boxes in Figure \ref{fig:boia})
and coordinate their effect.
More precisely, set $Z_t = \{s \in \BB^{k+1}\,|\,
(s,t) \not\in Y_2, f''(\tilde\Gamma_2(s,t)) \ne 0 \}$
for $t \in (0,2\pi)$.
The sets $Z_t$ form an open cover of the compact ball $\BB^{k+1}$.
Let $J = \{1, 2, \ldots, j_{\max}\}$
be a finite index set such that the sets $Z_{t_j}$ cover $\BB^{k+1}$.
Let $\breve K_j \subset K_j \subset Z_{t_j}$ be such that $K_j$ is compact
and the sets $\breve K_j$ form an open cover of $\BB^{k+1}$.
Let $\epsilon_0 > 0$ be such that the sets
$K_j \times [t_j - \epsilon_0, t_j + \epsilon_0]
\subset \BB^{k+1} \times \Ss^1$, $j \in J$,
are pairwise disjoint, disjoint from $Y_2$ and disjoint
from the set of pairs $(s,t)$ for which $f''(\tilde\Gamma_2(s,t)) = 0$.
Assume furthermore that $\epsilon_0 < \epsilon$,
where $\epsilon$ is given by Proposition \ref{prop:submersion}.
The sets $K_j \times T_j$, $T_j = [t_j - \epsilon_0, t_j + \epsilon_0]$,
are shown schematically in Figure \ref{fig:boia}.
In the figure, the dotted radii are the sets $\BB^{k+1} \times \{t_j\}$,
the two thick circles form the boundary $\Ss^k \times \Ss^1$
of $\BB^{k+1} \times \Ss^1$, the wiggly curve is $Y_2$
and the shaded boxes are $K_j \times T_j$:
notice that the interiors $\breve K_j$
indeed form an open cover in the figure.

We must first choose the bumps $\ell_i^j$: we will have
\[ \Gamma(s) = \Gamma_3(s) +
\sum_{i = 1, 2, 3;\; j \in J} a^j_i(s) \ell_i^j \]
where $a_i^j: \BB^{k+1} \to \RR$ are
continuous functions with support contained in $K_j$.
We next measure their capacity of correcting the monodromy.
Only then we choose $\breve Y_2$ so thin
that Lemma \ref{lemma:LinftyBeps} implies that
the monodromy of $\Gamma_3(s)$ is so close to $g$
that it can be fixed to obtain $\Gamma$ with the formula above.

Before proceeding we must clarify how the implicit function theorem
accounts for the capacity of a given triple of bumps to adjust monodromy.
As in Proposition \ref{prop:mufsurjective},
we use the $T$-variation $\mu_{f,T}$.
For $u \in L^\infty([0,2\pi])$ and $0 < t_- < t_+ < 2\pi$,
set $T_- = [0,t_-]$, $T_0 = [t_-,t_+]$, $T_+ = [t_+,2\pi]$,
$h_- = \mu_{f,T_-}(u)$, $h_0 = \mu_{f,T_0}(u)$ and $h_+ = \mu_{f,T_+}(u)$:
we have $\mu_f(u) = \mu_{f,[0,2\pi]}(u) = h_+h_0h_-$.
Now suppose that $\ell_i$, $i = 1,2,3$, are bumps supported in $T_0$
and consider $\tilde u = u + \sum_i a_i \ell_i$:
then $\mu_f(\tilde u) = h_+ \mu_{f,T_0}(\tilde u) h_-$.
We compare monodromies of $u$ and $\tilde u$ by writing
$\mu_{f,T_0}(\tilde u) = \psi \mu_{f,T_0}(u)$,
which is equivalent to
$\psi = (h_+)^{-1} \mu_f(\tilde u) (\mu_f(u))^{-1} h_+$.

\smallskip

{\nobf Claim:}
Let $K \subset \BB^{k+1}$ be a compact set
and $T = [t_-,t_+] \subset (0,2\pi)$.
Assume $K \times T$ to be disjoint from $Y_2$
and that $(s,t) \in K \times T$ implies $f''(\tilde\Gamma_2(s,t)) \ne 0$.
Assume also that $t_+ - t_- < \epsilon$, where $\epsilon$
is given in Proposition \ref{prop:submersion}.
Let $\ell_i$, $i = 1,2,3$, be bumps with disjoint supports contained in $T$
with $|\ell_i|_{L^\infty} < 1$.
Then there exists $\tilde\epsilon > 0$ such that,
for all continuous functions $\psi: K \to B_{\tilde\epsilon} \subset G$
with $\psi(s) = I$ for $s \in \partial K$
there are continuous functions $a_i: K \to [-E,E]$, $i = 1,2,3$,
with $a_i(s) = 0$ for $s \in \partial K$ for which
\[ \mu_{f,T}\left(\Gamma_2(s) + \sum_{i = 1,2,3} a_i(s) \ell_i\right)
= \psi(s) \mu_{f,T}(\Gamma_2(s)). \]

{\nobf Proof of claim:}
Define $\zeta: K \times \RR^3 \to G$ by
\[ \zeta(s,a_1,a_2,a_3) =
 \mu_{f,T}\left(\Gamma_2(s) + \sum_i a_i \ell_i\right)
(\mu_{f,T}(\Gamma_2(s)))^{-1}.  \]
Clearly $\zeta(s,0,0,0) = I$ for all $s \in K$.
From Proposition \ref{prop:submersion}, the three vectors
\( {\partial \zeta_j}/{\partial a_i} \)
are linearly independent.
From the implicit function theorem, there exists
$\tilde\epsilon > 0$ and a continuous function
$\alpha: K \times B_{\tilde\epsilon} \to \RR^3$
with $\alpha(s,I) = 0$ and
$\zeta(s,\alpha(s,h)) = h$ for all $s \in K_j$
and $h \in B_{\tilde\epsilon}$.
We can furthermore assume without loss of generality
that $||\alpha(s,h)|| < E$ for all $s \in K$
and $h \in B_{\tilde\epsilon}$.
Now take $(a_1(s), a_2(s), a_3(s)) = \alpha(s,\psi(s))$,
completing the proof of the claim.
\qede

\smallskip

We now put the partial results to work:
for each $j \in J$, construct bumps $\ell^j_1, \ell^j_2, \ell^j_3$
with support contained in $T_j$ and $|\ell^j_i|_{L^\infty} < 1$.
Apply the claim to $K_j$, $T_j$ and $\ell^j_i$ to obtain $\tilde\epsilon_j$;
let $\epsilon_1$ be the minimum among all $\tilde\epsilon_{j}$, $j \in J$.
The set of all $\tilde\Phi(t)$ where
\[ \tilde\Phi(0) = I, \quad
\tilde\Phi'(t) =
\begin{pmatrix} 0 & 1 \\ f'(u(t)) & 0 \end{pmatrix} \tilde\Phi(t), \quad
t \in [0,2\pi], \quad |u|_{L^\infty(\Ss^1)} < 2E, \]
is contained in a compact set $K' \subset G$.
Let $\epsilon'_2 > 0$ be such that
\[ B_{\epsilon'_2} \subseteq \bigcap_{h \in K'} h^{-1} B_{\epsilon_1} h \]
and $\epsilon_2 > 0$ such that
$B_{\epsilon_2} B_{\epsilon_2} \subseteq B_{\epsilon'_2}$.
From Lemma \ref{lemma:LinftyBeps}, let $\epsilon_3 > 0$ be such that
\[ |u_0|_{L^\infty}, |u_1|_{L^\infty} < 2E, \quad
\lambda(\{t | u_0(t) \ne u_1(t)\}) < \epsilon_3 \quad\Rightarrow\quad
\mu_f(u_0)(\mu_f(u_1))^{-1} \in B_{\epsilon_2}. \]

Select an open set $\breve Y_2 \subset \BB^{k+1} \times \Ss^1$,
$Y_2 \subset \breve Y_2$,
which is removed from the boundary and satisfies
$\lambda(\breve Y_2(s)) < \epsilon_3$ for all $s \in \BB^{k+1}$ where,
as before,
\[\breve Y_2(s) =
\Pi_{\Ss^1}\left(\breve Y_2 \cap (\{s\} \times \Ss^1)\right). \]
Define a continuous function $\tilde\Gamma_3$
with $|\tilde\Gamma_3(s,t)| < 2E$ for all $s, t$
and a corresponding $\Gamma_3: \BB^{k+1} \to C^0(\Ss^1)$
coinciding with $\tilde\Gamma_2$ outside $\breve Y_2$
(and otherwise arbitrary in $\breve Y_2$).
From the construction of $\epsilon_3$, we have
$\mu_f(\Gamma_3(s)) g^{-1} \in B_{\epsilon_2}$ for all $s$
(recall that $\mu_f(\Gamma_2(s)) = g$).

We will define inductively in $j$ the functions $a_i^j$ and
\[ \Gamma_{4,j}(s) = \Gamma_{4,j-1}(s) +
\sum_{i = 1, 2, 3} a^{j}_i(s) \ell_i^{j}, \quad
\Gamma_{4,0} = \Gamma_3; \]
notice that $\Gamma_{4,j_{\max}} = \Gamma$.
The idea is that we activate one $K_j \times T_j$ box at a time.

Recall that $\expg: \fg \to G$ is the exponential map.
Let $v_0: \BB^{k+1} \to \fg$ be defined by
$\mu_f(\Gamma_3(s)) g^{-1} = \expg(v_0(s))$;
the function $v_0$ provides a linear measure for
the error in the monodromy of $\Gamma_3$:
we now construct intermediate functions $v_j: \BB^{k+1} \to \fg$
with $v_{j_{\max}} = 0$ and functions $\Gamma_{4,j}$ with
$\mu_f(\Gamma_{4,j}(s)) g^{-1} = \expg(v_j(s))$.
Let $r_j: \BB^{k+1} \to [0,1]$ be a smooth partition of unity
associated to the partition $\breve K_j$, so that $\sum_j r_j = 1$
and the support of $r_j$ is contained in $\breve K_j$.
Set
\[ v_j(s) = \left( \sum_{j' > j} r_{j'}(s) \right) v_0(s); \]
notice that $v_j(s) = v_{j-1}(s)$ for $s \notin \breve K_j$.
Define inductively in $j$ the functions $a_i^j$ so that
\[ \Gamma_{4,j}(s) = \Gamma_{4,j-1}(s) +
\sum_{i = 1, 2, 3} a^{j}_i(s) \ell_i^{j}, \quad
\Gamma_{4,0} = \Gamma_3, \quad
\mu_f(\Gamma_{4,j}(s)) g^{-1} = \expg(v_j(s)): \]
then set $\Gamma_{4,j_{\max}} = \Gamma$.
We are left with showing that this is indeed possible,
i.e., that the functions $a_i^j$ are continuous and well defined.

Assume that $\Gamma_{4,j-1}$ has been constructed,
in other words, the functions $a^{j'}_i$ have been obtained for $j' < j$
and we have $\mu_f(\Gamma_{4,j - 1}(s)) g^{-1} = \expg(v_{j-1}(s))$.
We need functions $a^j_i$ such that
\[ \mu_f\left(\Gamma_{4,j - 1}(s) +
\sum_{i = 1,2,3} a^j_i(s) \ell^j_i \right) g^{-1} = \expg(v_j(s)). \]
Notice that $a^j_i(s) = 0$ for $s \not\in K_j$, as required.
For $s \in K_j$ and $T_j = [t_j^-, t_j^+]$, set
$h_-(s) = \mu_{f,[0,t_j^-]}(\Gamma_{4,j-1}(s))$,
$h_0(s) = \mu_{f,T_j}(\Gamma_{4,j-1}(s))$,
$h_+(s) = \mu_{f,[t_j^+,2\pi]}(\Gamma_{4,j-1}(s))$
so that
\[ h_+(s) h_0(s) h_-(s) = \mu_f(\Gamma_{4,j-1}) = \expg(v_{j-1}(s)) g. \]
Define $\psi_j: \BB^{k+1} \to G$ so that
$h_+(s) \psi_j(s) h_0(s) h_-(s) = \expg(v_j(s)) g$;
in other words,
\[ \psi_j(s) = (h_+(s))^{-1} \expg(v_j(s)) \expg(-v_{j-1}(s)) h_+(s). \]
By construction, $h_+(s) \in K' \subset G$ and
$\expg(v_j(s)), \expg(-v_{j-1}(s)) \in B_{\epsilon_2}$ whence
$\expg(v_j(s))\expg(-v_{j-1}(s)) \in B_{\epsilon'_2}$,
so that $\psi_j(s) \in B_{\epsilon_1}$.
Apply the claim to $K_j$, $T_j$, $\ell^j_i$ and $\psi_j$
to obtain $a^j_i$ and we are done.
\qed

\section{Hilbert manifolds with cone-like singularities}

In this section we obtain normal forms near the critical set $\Cc$.
We first consider a simple scenario.
Recall that $\cone\; \subset \RR^3$ is the cone $x^2 + y^2 - z^2 = 0$
with vertex $0 \in \RR^3$.

\begin{theo}
\label{theo:danA}
Let $\Hh$ be an infinite dimensional smooth Hilbert manifold
and $\sigma: \Hh \to \RR^3$ a smooth surjective submersion
with contractible fibers.
Set $\Cc_1 = \sigma^{-1}(\cone)$ and $\Cc_2 = \sigma^{-1}(\{0\})$.
Then there exists a diffeomorphism
\[ \theta: (\Hh, \Cc_1, \Cc_2) \to (\RR^3, \cone, \{0\}) \times \HH. \]
\end{theo}

The proof will be organized in two steps
which are stated now and justified later.
For a manifold with boundary $X$ (of finite or infinite dimension),
we write $\interior X$ for its interior and $\partial X$ for its boundary.

\noindent {\bf Step 1} \textit{There is a closed tubular neighborhood
$\Dd \subset \Hh$ of the smooth submanifold $\Cc_2$ and a diffeomorphism
$\theta_1: (\Dd, \partial \Dd) \to (D_1, \partial D_1) \times \Cc_2$ satisfying
$\theta_1(\Cc_2) = 0 \times \Cc_2$ and
$\theta_1(\Dd \cap \Cc_1) = (D_1 \cap \cone) \times \Cc_2$.  }

Here $D_r$ denotes the closed disk (or ball) of radius $r$ in $\RR^3$
centered at the origin.
Recall that a tubular neighborhood $\Dd$ is the image by an embedding
of the \textit{unit disk bundle},
a subset of the normal bundle of $\Cc_2$ in $\Hh$, $\nu: \Dd \to \Cc_2$,
which, in this case, is trivial since $\Cc_2$ is contractible.
The fibers of $\nu: \Dd \to \Cc_2$ can be identified to the unit disk in
$\RR^3$, and therefore $\nu$ can be identified to
$\Pi_1: \Cc_2 \times D_1 \to \Cc_2$,
where $\Pi_1$ is the projection on the first coordinate.
Note that $\Dd$ is a (codimension zero) smooth submanifold with boundary,
$\partial \Dd$ being diffeomorphic to $\Cc_2 \times \Ss^2$.

Since $\sigma$ is a submersion and, for $\Dd$ satisfying Step 1,
$\cone$ is transversal to $\partial D_1$,
we have that $\Cc_1$ is also transversal to $\partial \Dd$.
We then consider the smooth Hilbert manifold with boundary
$\Vv = \Hh \smallsetminus \interior \Dd$.
The subset $\Kk = \Cc_1 \smallsetminus  \interior \Dd$ is a
codimension $1$ submanifold with boundary
which intersects $\partial \Vv = \partial \Dd$ transversally.

\noindent {\bf Step 2} \textit{There is a diffeomorphism
$\theta_2: (\Vv, \partial \Vv) \to
(\RR^3 \smallsetminus \interior D_1, \partial D_1) \times \Cc_2$ with
$\theta_2(\Kk) = (\cone \smallsetminus \interior D_1) \times \Cc_2$.}

We may furthermore assume that $\theta_1$ and $\theta_2$,
given by Steps 1 and 2, coincide on $\partial\Vv = \partial\Dd$.
Indeed, start with $\theta_1$ given by Step 1;
in order to define $\theta_2$ we first set
$\theta_2|_{\partial\Vv} = \theta_1|_{\partial\Vv}$,
$\theta_2|_{\partial\Vv}: (\partial \Vv, \partial \Vv \cap \Cc_1) \to
(\partial D_1, \partial D_1 \cap \cone) \times \Cc_2$.

From Step 2 and the identification between
$(\RR^3 \smallsetminus  \interior D_1,
\cone \smallsetminus  \interior D_1)$ and
$(\partial D_1, \cone \cap\;\partial D_1) \times [1,\infty)$,
extend this restriction to obtain the desired diffeomorphism $\theta_2$.
Combining such $\theta_1$ and $\theta_2$,
and in view of the fact that $\Cc_2$ is diffeomorphic to $\HH$,
one obtains Theorem \ref{theo:danA}.

We need a few additional notations and observations.
For a positive smooth function $\beps: \Cc_2 \to \RR_+$,
we denote by
$(\Cc_2 \times D)_\beps = \{ (p,v)\;|\;||v|| \le \beps(p) \}$.
Let $\Pi_1: (\Cc_2 \times D)_\beps \to \Cc_2$
be the first factor projection,
whose fiber above $p \in \Cc_2$ is $p \times D_{\beps(p)}$.
Note that if $\beps$ is a smooth function then
$(\Cc_2 \times D)_\beps$ is a smooth Hilbert manifold with boundary
and one can produce a fiber diffeomorphism
$\theta_\beps: (\Cc_2 \times D)_\beps \to (\Cc_2 \times D)_\bone$
where $\bone$ denotes the constant function and
the diffeomorphism keeps the first coordinate fixed.

To accomplish Step 1, we use a smooth closed tubular neighborhood
$\omega: (\Cc_2 \times D)_\beps \to \Hh$,
where $\beps: \Cc_2 \to \RR_+$ is a smooth function
so that the composition $\sigma \circ \omega$,
when restricted to the fibers $p \times D_{\beps(p)}$
is a diffeomorphism onto its image.
Clearly, for such $\omega$,
$(\sigma \circ \omega, \Pi_1): (\Cc_2 \times D)_\beps \to \RR^3 \times \Cc_2$
is a diffeomorphism onto the image.
To construct such pair $(\omega, \beps)$, we choose
$\omega_0 : \Cc_2 \times \RR^3 \to T(\Hh)|_{\Cc_2}$, a splitting of the
surjective bundle map
$D\sigma: T(\Hh)|_{\Cc_2} \to \Cc_2^\ast(T(\RR^3))|_{\Cc_2}
= \Cc_2 \times \RR^3$,
fixing base point, for which $D\sigma \circ \omega_0 = \id$. Such
splitting exists since $\sigma$ is a submersion.
We choose a complete
Riemannian metric on $\Hh$ and use the exponential map with respect to
this metric to define the smooth map $e: \Cc_2 \times \RR^3 \to \Hh$
(we regard $\Cc_2 \times \RR^3$ as the total space of the normal bundle of
$\Cc_2$ in $\Hh$). The differential of $e$ at any point of $\Cc_2 \times 0$
(the zero section of the trivial bundle  $\Cc_2 \times \RR^3 \to \Cc_2)$ is
an isomorphism and therefore there exists $\beps_1: \Cc_2 \to \RR_+$
so that $e$ restricted to $(\Cc_2 \times D)(\beps_1)$ is a
diffeomorphism onto its image. Since the differential at zero of
$\sigma \circ e$ when restricted to any fiber of $p \times \RR^3 \to \Cc_2$ is
an isomorphism, one obtains $\beps_2: \Cc_2 \to \RR_+$ so that
$\sigma \circ e$, when restricted to the disk $D_{\beps_2(p)}$, is a
diffeomorphism on the image.
Choose $\beps \leq \min(\beps_1, \beps_2)$
but still smooth and positive, and take $\omega$ the
restriction of $e$ to $(\Cc_2 \times D)_\beps$,
thus completing the proof of Step 1.

For Step 2, we need a few preliminary results.
We begin with an easy consequence of item 1 of Proposition 3.1 in \cite{BST1}.

\begin{lemma}
\label{lemma:dan62empty}
Suppose $(\Vv, \partial \Vv)$ is a Hilbert manifold with boundary
and that the inclusion $\partial \Vv \hookrightarrow \Vv$
is a homotopy equivalence.
Then there exists a diffeomorphism $\theta:(\Vv, \partial \Vv) \to
(\partial \Vv \times [1, \infty), \partial \Vv \times \{1\})$.
\end{lemma}

We now use this lemma to prove an amplification
(the lemma is the degenerate case $\Kk = \emptyset$).

\begin{prop}
\label{prop:dan62}
Suppose $(\Vv, \partial \Vv)$ is a Hilbert manifold with boundary, 
$\Kk \subset \Vv$ is a finite codimension submanifold such that $\Kk$ is
transversal to $\partial \Vv$ and $\partial \Kk = \Kk \cup \partial \Vv$.
Suppose also that the following inclusions are homotopy equivalences:
$\partial \Vv \hookrightarrow \Vv$, $\partial \Kk \hookrightarrow \Kk$,
$\partial \Vv \smallsetminus  \partial \Kk \hookrightarrow \Vv \smallsetminus  \Kk$.
Then there exists a diffeomorphism 
\( \theta:(\Vv, \partial \Vv) \to
(\partial \Vv \times [1, \infty), \partial \Vv \times \{1\}) \)
so that $\theta(\Kk) = \partial \Kk \times [1, \infty)$.
\end{prop}

{\nobf Proof: }
Construct a relative collar neighborhood $\phi$
of $(\partial \Vv, \partial \Kk)$ in $(\Vv,\Kk)$, i.e.,
a closed embedding
$\phi: ( \partial \Vv \times [1,2], \partial \Vv \times \{1\}) \to
(\Vv, \partial \Vv)$ such that its restriction
to $\partial \Vv \times \{1\}$ is the projection on the first coordinate and
$\phi(\partial \Vv \times [1,2]) \cap \Kk = \phi( \partial \Kk \times [1,2])$.
It follows from lemma \ref{lemma:dan62empty}
and from the homotopy equivalence $\partial \Kk \hookrightarrow \Kk$ that
the pair $( \partial \Kk \times [1, \infty), \partial \Kk \times \{1\})$
is diffeomorphic to $(\Kk, \partial \Kk)$:
we may assume that this diffeomorphism coincides with $\phi$
on $\partial \Kk \times [1,2]$ and we also call it $\phi$.
Let $\Kk^1 = \phi(\partial \Kk \times [1,2))$ and
$\Kk^2 = \phi(\partial \Kk \times [2, \infty))$.
Standard arguments (``pushing to infinity") imply that the triple
$(\Vv, \partial \Vv, \Kk)$ is diffeomorphic to
$(\Vv \smallsetminus  \Kk^2, \partial \Vv, \Kk^1)$.
The set $\Ww = (\Vv \smallsetminus  \Kk^2) \smallsetminus
\phi(\partial \Vv \times [1,2))$ is a smooth manifold with boundary,
$\partial W$ being
diffeomorphic to $( \partial \Vv \smallsetminus  \Kk) \times \{2\}$.
The homotopy equivalences in the statement
imply that $\partial \Ww \to \Ww$ is a homotopy equivalence and then,
by lemma \ref{lemma:dan62empty}, we have that $(\Ww, \partial \Ww)$ and
$(\partial \Ww \times [2, \infty), \partial \Ww \times \{2\})$ are
diffeomorphic. We conclude that there exists a diffeomorphism from
$( \Vv \smallsetminus  \Kk^2, \partial \Vv)$ to
$( \partial \Vv \times [1, \infty) \smallsetminus
\partial \Kk \times [2,\infty), \partial \Vv \times \{1\})$ which sends
$\Kk^1$ into $\partial \Kk \times [1,2)$.
By the same trick of pushing to infinity,
there is a diffeomorphism from the last pair to
$(\partial \Vv \times [1, \infty), \partial \Vv \times \{1\})$
which sends $\partial \Kk \times [1,2)$
into $\partial \Kk \times [1, \infty)$, concluding the proof.
\qed

Step 2 is now accomplished
and so is the proof of Theorem \ref{theo:danA}.
In order to prove Theorem \ref{theo:main} in the next section,
we need to strengthen Theorem \ref{theo:danA} somewhat:
this is done in Propositions \ref{prop:danB} and \ref{prop:dan66} below.
A \textit{closed disk} $\Bb$ in a Hilbert manifold is a closed tubular
neighborhood of a point, i.e.,
a closed set which is a codimension zero submanifold
with contractible boundary.

\begin{prop}
\label{prop:danB}
Let $\Hh$ be an infinite dimensional smooth Hilbert manifold
and $\sigma: \Hh \to \RR^3$ a smooth surjective submersion
with contractible fibers.
Let $\Bb_1$ and $\Bb_2$ (resp. $\Bb_1'$ and $\Bb_2'$)
be two closed disks contained
in the two contractible components of $\Hh \smallsetminus  \Cc_1$
(resp.  $(\RR^3 \smallsetminus  \cone) \times \HH$).
Then there exists a diffeomorphism
$\theta: (\Hh, \Cc_1, \Cc_2) \to (\RR^3, \cone, 0) \times \HH$
so that $\theta(\Bb_i) = \Bb_i'$, $i=1,2$.
\end{prop}

{\nobf Proof: }
First construct
$\theta_1: (\Hh, \Cc_1, \Cc_2) \to (\RR^3, \cone, 0) \times \HH$
using Theorem \ref{theo:danA}.
The diffeomorphism $\theta_1$ sends $\Cc_1$ into $\cone \times \HH$ and
can be easily modified away from a neighborhood of $\Cc_1$ to obtain
the desired diffeomorphism $\theta$ in view of Fact \ref{fact:cerf} below.
\qed

\begin{fact}[Cerf lemma]
\label{fact:cerf}
Let $\Bb$ be the closed disk in $\HH$,
$\phi_1, \phi_2: \Bb \to \Hh$ be smooth codimension zero closed embeddings
and $\cal{U} \subset \Hh$ an open set with $\phi_i(\Bb) \subset \cal{U}$.
Then there is an isotopy $h_t: \Hh \to \Hh$, $t \in [0,1]$,
with $h_t|_{\Hh \smallsetminus  \cal{U}} = \id$,
$h_0 = \id$ and $h_1 \circ \phi_2 = \phi_1$.
\end{fact}

See \cite{Cerf} and \cite{Palais}
for the finite dimensional case of Cerf lemma
(or \cite{Hirsch}, Theorem 3.1 in Chapter 8);
see \cite{BK} for the (similar) infinite dimensional case.

The following lemma is another easy consequence of
Proposition 3.1 in \cite{BST1}.

\begin{lemma}
\label{lemma:dan65}
Let $\Vv_0$ be a contractible Hilbert manifold.
Let $\Vv_1 = \Vv_0 \times [0,1]$,
$\Vv_2 = \Vv_0 \times [0,1) \smallsetminus  \interior\Bb$,
where $\Bb$ is a closed disk contained in $\Vv_0 \times (0,1)$.
Also, $ \partial \Vv_1 = \Vv_0 \times \{ 0,1\}$ and
$\partial \Vv_2 = \Vv_0 \times \{0\} \cup \partial \Bb$.
Then there exists a diffeomorphism
$\theta: (\Vv_1, \partial \Vv_1) \to (\Vv_2, \partial \Vv_2)$
such that $\theta|_{\Vv_0 \times \{0\}} = \id$.
\end{lemma}

\begin{prop}
\label{prop:dan66}
Let $(\Vv, \partial \Vv)$ be a Hilbert manifold with boundary consisting
of two contractible components $\partial_- \Vv$ and $\partial_+ \Vv$.
Let $\phi_{\pm}: \partial_{\pm} \Vv \times [0,1] \to \Vv$ be two
disjoint closed  collar neighborhoods of $\partial_\pm \Vv$, say
$\Dd_{\pm} = \phi_{\pm}(\partial_{\pm} \Vv \times [0,1])$
and let $\Bb_\pm$ be two closed disks contained in
$\phi_{\pm} ( \partial_{\pm} \Vv \times (0,1))$.
Then there exists a diffeomorphism
$\theta: (\Vv, \partial \Vv) \to
(\Vv \smallsetminus (\partial \Vv \cup \interior(\Bb_+ \cup \Bb_-)),
\partial (\Bb_- \cup \Bb_+))$
which restricts to the identity on $\Vv \smallsetminus (\Dd_+ \cup \Dd_-)$.
\end{prop}

\noindent\textbf{Proof:}
Apply lemma \ref{lemma:dan65} twice: in the first instance, take
$\Vv_1 = \phi_+(\partial_+ \Vv \times [0,1])$,
$\partial\Vv_1 = \phi_+(\partial_+ \Vv \times \{0,1\})$,
$\Vv_2 = \Vv_1 \smallsetminus
(\phi_+(\partial_+ \Vv \times \{0\}) \cup \interior \Bb_+)$,
$\partial\Vv_2 = \phi_+(\partial_+ \Vv \times \{1\}) \cup \partial\Bb_+$.
In the second, replace the $+$ signs by $-$.
\qed

\goodbreak

\section{Proof of the main theorem}

The idea of the proof of Theorem \ref{theo:main}
is to consider the restriction of $\mu_f$ to the preimage
of the set of lifted matrices with nonnegative trace
(the regions between dotted vertical lines containing components
of $T_2$ in Figure \ref{fig:g0p})
and use Proposition \ref{prop:dan66}.
Extending the diffeomorphism to the complement is then easy.
To insure that the hypothesis of Proposition \ref{prop:dan66} hold,
however, we need Michael's continuous selection theorem
(the  main result in \cite{Michael}), which we state below
in the special situation we use.

Recall that a surjective continuous map $f:X\to Y$ of arbitrary
topological spaces is a \textit{topological submersion}
if for each $x\in X$ there exist
a neighborhood $U\subset Y$ of $f(x)$,
a neighborhood $V\subset f^{-1}(f(x))$ of $x$ and
an open embedding $h:U\times V\to X$ so that
$f\circ h$ is the projection on the first coordinate.
Clearly a smooth surjective submersion of Hilbert manifolds is
a topological submersion and the  pullback of
a topological submersion by a continuous map is a topological submersion.

\begin{fact}[Michael's theorem]
\label{fact:michael}
Let $X$ and $Y$ be (possibly infinite dimensional) manifolds
and $f:X\to Y$  be a topological submersion with $k$-connected fibers.
Then $f$ induces an isomorphism on homotopy groups
in dimension smaller than $k$ and an epimorphism in dimension $k$.
In particular, if the fibers are contractible
then $f$ is a homotopy equivalence.
\end{fact}

The proof of this result is based on the observation that for a 
topological submersion $f$ the  assignment $y\in Y \rightsquigarrow f^{-1}(y)$
is lower semi-continuous in the sense of Michael (since $f$ is open) and
equi-$LC^m$ (locally $m$-connected) for any $m$
(see \cite{Michael} for definitions).
In view of the main theorem in \cite{Michael}, if $Y$ has dimension $n$ and
any fiber is $n$-connected then any continuous section
$\varphi : A\to X$, where $A$ is a closed subset of $Y$,
has an extension to a section $\tilde\varphi:Y\to X$.
This clearly implies the first part of the statement
and therefore the second, since manifolds have the homotopy type of ANR's.

The proposition below gives a normal form for a submersion
with contractible fibers onto a closed neighborhood of a cone.

\begin{prop}
\label{prop:dan64}
Let $n \in \ZZ$, $n > 0$, set $I_n = [2\pi n - \pi/2,2\pi n + \pi/2]$
and $\cone_n \; =
\{ (x,y,z) \in \RR^2 \times \interior I_n \,|\,
x^2 + y^2 = \tan^2 z \}$.
Let $(\Hh, \partial \Hh)$ be a Hilbert manifold with boundary
and $\sigma : \Hh \to \RR^2 \times I_n$
be a surjective submersion with contractible fibers so that
$\partial \Hh = \sigma^{-1}(\RR^2 \times \partial I_n)$.
Let $\Cc = \sigma^{-1}(\cone_n)$.
Then there exists a diffeomorphism $\theta: (\Hh, \partial \Hh) \to
\RR^2 \times (I_n, \partial I_n) \times \HH$
so that $\theta(\Cc) = \;\cone_n \times \HH$.
\end{prop}

\noindent\textbf{Proof:}
Set $\Vv = \Hh$,
$\partial_{\pm} \Vv = \sigma^{-1}(\RR^2 \times \{2\pi n \pm \pi/2\})$.
Fact \ref{fact:michael} implies that each space
$\partial_{\pm} \Vv$ is contractible.
Choose two disjoint closed collar
neighborhoods $\Dd_{\pm}$ of $\partial_{\pm} \Vv$ which are
disjoint from $\Cc$ and two disjoint closed collar neighborhoods
$\Dd_{\pm}'$ of $\RR^2 \times \{ 2\pi n \pm \pi/2\} \times \HH$
which are disjoint from $\cone_n \times \HH$.
Now, choose closed disks $\Bb_+$, $\Bb_-$, $\Bb_+'$ and $\Bb_-'$
in the interior of the closed tubular neighborhoods
$\Dd_+$, $\Dd_-$, $\Dd_+'$ and $\Dd_-'$, respectively.
Set $\Bb_\pm = \Bb_- \cup \Bb_+$ and
$\Bb_\pm' = \Bb_-' \cup \Bb_+'$.
Use Proposition \ref{prop:dan66} to construct the diffeomorphisms
\begin{gather*}
\theta_1: (\Vv, \partial \Vv) \to
(\Vv \smallsetminus ( \partial \Vv \cup \interior(\Bb_\pm)),
\partial (\Bb_\pm)), \\
\theta_2: ( I \times \RR^2 \times \HH,
\partial I \times \RR^2 \times \HH) \to
( (\interior I \times \RR^2 \times \HH) \smallsetminus
(\interior(\Bb_\pm')),
\partial (\Bb_\pm')) 
\end{gather*}
and Proposition \ref{prop:danB} to construct the diffeomorphism
\[ \theta:(\Vv \smallsetminus
(\partial \Vv \cup \interior(\Bb_\pm)),
\partial(\Bb_\pm)) \to
(\interior I \times \RR^2 \times \HH \smallsetminus
\interior(\Bb_\pm'),
\partial (\Bb_\pm')). \]
The desired diffeomorphism 
is $\theta_2^{-1} \circ \theta \circ \theta_1$.
\qed

We now return to the nonlinear monodromy $\mu_f: X^\ast \to G_0$
where $X = H^p(\Ss^1)$.
From Proposition \ref{prop:submersion}, $\mu_f$ is a submersion
provided that $f'$ is nowhere flat.
If additionally $f'$ is surjective,
Proposition \ref{prop:mufsurjective} implies the surjectivity of $\mu_f$ and
Theorem \ref{theo:contract} proves that
the fibers of $\mu_f$ are contractible.
Let $Z_2 = \Pi^{-1}(I) \subset T_2 \subset G$ be the set of vertices
of cones in $T_2$: notice that $Z_2$ is a group isomorphic to $\ZZ$.

\begin{prop}
\label{prop:main1}
Let $f: \RR \to \RR$ be a smooth nonlinearity.
Assume that $f'$ is surjective and nowhere flat.
Let $X = H^p(\Ss^1)$, $p \ge 1$,
$X^\ast = X \smallsetminus \{ \textrm{constant functions} \}$
and $\mu_f: X^\ast \to G_0$ be the monodromy map.
Let $\Cc^\ast = \Cc \cap X^\ast = \mu_f^{-1}(G_0 \cap T_2)$
and $\Cc^\ast_2 = \mu_f^{-1}(G_0 \cap Z_2)$.
Then the triple $(X^\ast, \Cc^\ast, \Cc^\ast_2)$
is diffeomorphic to the triple
$(G_0, G_0 \cap T_2, G_0 \cap Z_2) \times \HH$.
\end{prop}

\noindent\textbf{Proof:}
Let $\Sigma \subset \RR^3$ be as defined in the introduction
and $\Sigma_2 = \{ (0,0,2\pi n), n \in \ZZ, n > 0\} \subset \Sigma$
the set of vertices of cones in $\Sigma$.
From \cite{BST2} there exists a diffeomorphism
$\psi: (G_0, G_0 \cap T_2, G_0 \cap Z_2)
\to (\RR^3, \Sigma, \Sigma_2)$.
Apply Proposition \ref{prop:dan64} for $\sigma = \psi \circ \mu_f$,
$\Hh = \sigma^{-1}(\RR^2 \times I_n)$;
attaching the pieces presents no difficulty.
\qed

The attentive reader will notice that this proposition also holds
(with the same proof) if $X$ is a smoothing Hilbert space
or if it is a separable Hilbert space satisfying the hypothesis
of corollary \ref{coro:submersionX}.
Also, if $X$ is a smoothing Banach space
(or a Banach space satisfying the hypothesis
of corollary \ref{coro:submersionX})
then a similar but weaker conclusion holds:
the triples are homeomorphic.

The use of Michael's Theorem can be avoided if we are willing
to prove directly that if $P \subset G_0$ is diffeomorphic to a plane then
the sets $\mu_f^{-1}(P) \cap X^\ast$ are contractible.
This would be achieved by mimicking the proof 
of Theorem \ref{theo:contract}.

Theorem \ref{theo:main} is slightly different from
Proposition \ref{prop:main1}: the Proposition is stated
for the triple $(X^\ast, \Cc^\ast, \Cc^\ast_2)$ while
the theorem is stated for the pair $(X, \Cc^\ast)$.
Proposition \ref{prop:localmu} guarantees that if $f$ is good
then $\Cc^\ast = \Cc$.
Also, $f$ is admissible 
then $\Cc$ is the disjoint union of $\Cc^\ast$ and isolated points.

\bigbreak


%% file: periodicX.tex

\bigbreak

{

\parindent=0pt
\parskip=0pt
\obeylines

Dan Burghelea, Ohio State University, burghele@math.ohio-state.edu
Nicolau C. Saldanha, PUC-Rio, nicolau@mat.puc-rio.br
Carlos Tomei, PUC-Rio, tomei@mat.puc-rio.br

\bigskip

Department of Mathematics, Ohio State University,
231 West 18th Ave, Columbus, OH 43210-1174, USA

\smallskip

Departamento de Matem\'atica, PUC-Rio
R. Marqu\^es de S. Vicente 225, Rio de Janeiro, RJ 22453-900, Brazil

\smallskip


}

%% file: periodic.bbl
\bigskip

\begin{thebibliography}{[10]}

\bibitem{AM}{ Ambrosetti, A. and Malchiodi, A.,
{\sl Nonlinear Analysis and Semilinear Elliptic Problems},
Cambridge studies in advanced mathematics, 104,
CUP, Cambridge, UK, 2007.}
\bibitem{AP}{ Ambrosetti, A. and Prodi, G.,
{\sl On the inversion of some differentiable maps between Banach spaces
with singularities},
Ann. Mat. Pura Appl. 93, 231-246, 1972.}
\bibitem{BC} { Berger, M.S. and Church, P.T., {\sl Complete integrability
and perturbation of a nonlinear Dirichlet problem},
Ind. Univ. Math. J. 28, 935-952, 1979.}
\bibitem{BP}{ Berger, M. and Podolak, E.,
{\sl On the solutions of a nonlinear Dirichlet problem},
Ind. Univ. Math. J. 24, 837-846, 1975.}
\bibitem{BT}{ Bueno, H. and Tomei, C.,
{\sl Critical sets of nonlinear Sturm-Liouville operators
of Ambrosetti-Prodi type},
Nonlinearity 15, 1073-1077, 2002.}
\bibitem{BK}{ Burghelea, D. and Kuiper, N.,
{\sl Hilbert manifolds},
Ann. of Math. 90, 379-417, 1969.}
\bibitem{BST1}{ Burghelea, D., Saldanha, N. and Tomei, C.,
{\sl Results on infinite dimensional topology and
applications to the structure of the critical set
of non-linear Sturm-Liouville operators},
J. Differential Equations, 188, 569-590, 2003.}
\bibitem{BST2}{ Burghelea, D., Saldanha, N. and Tomei, C.,
{\sl The topology of the monodromy map of the second order ODE},
J. Differential Equations, 227, 581-597, 2006.}
\bibitem{Cerf}{Cerf, J.,
{\sl Topologie de certains espaces de plongements},
Bull. Soc. Math. France, 89, 227-380, 1961.}
\bibitem{CT}{ Church, P. T. and Timourian, J. G.,
{\sl Global structure for nonlinear
operators in differential and integral equations. I. Folds, II. Cusps},
Topological nonlinear analysis, II (Frascati, 1995), 109--160, 161--245,
Progr. Nonlinear Differential Equations Appl., 27, Birkh\"auser,
Boston, MA, 1997.}
\bibitem{Hirsch}{ Hirsch, M. W.,
{\sl Differential Topology},
GTM 33, Springer, 1976.}
\bibitem{MST0}{ Malta, I., Saldanha, N. C. and Tomei, C.,
{\sl The numerical inversion of functions from the plane to the plane},
Math. Comp. 65, no. 216, 1531-1552, 1996.}
\bibitem{MST2}{ Malta, I., Saldanha, N. C. and Tomei, C.,
{\sl Regular level sets of averages of
Nemytski{\u\i} operators are contractible},
J. Func. Anal., 143, 461-469, 1997.}
\bibitem{MST1}{ Malta, I., Saldanha, N. C. and Tomei, C.,
{\sl Morin singularities and global geometry in a class
of ordinary differential equations},
Top. Meth. in Nonlinear Anal., 10 (1), 137-169, 1997.}
\bibitem{Michael}{ Michael, E. A.,
{\sl Continuous selections III},
Annals of Mathematics, 65 (2), 375-390, 1957.}
\bibitem{Palais}{ Palais, R.,
{\sl Extending diffeomorphisms},
Proc. Amer. Math. Soc. 11, 274-277, 1960.}
\bibitem{PT}{ P\"oschel, J. and Trubowitz, E.,
{\sl Inverse spectral theory},
Academic Press, Boston, 1987.}
\bibitem{Ruf}{ Ruf, B.,
{\sl Singularity theory and bifurcation phenomena
in differential equations},
quaderno 19, Univ. Studi Milano, 1996.}
\end{thebibliography}
